# Potential flow in the critical strip and the Riemann hypothesis

J. G. Andrews

H.H. Wills Physics Laboratory, University of Bristol, Bristol, United Kingdom
john.andrews@bristol.ac.uk

ABSTRACT. We describe the behaviour of the Dirichlet eta function in the critical strip, in terms of the potential flow of an ideal fluid. Using well-known results from complex potential theory and number theory, we show that the Dirichlet eta function has no zeros in the critical strip off the critical line, consistent with the Riemann hypothesis.

## 1    Introduction

A function of a complex variable

$$f(\xi) = \phi(x, y) + i\psi(x, y) \qquad \xi = x + iy \tag{1}$$

is analytic provided its derivative $f'(\xi)$ exists, and this requires that $\phi(x,y)$ and $\psi(x,y)$ satisfy the Cauchy-Riemann conditions

$$\frac{\partial \phi}{\partial x} = \frac{\partial \psi}{\partial y} \; , \quad \frac{\partial \phi}{\partial y} = -\frac{\partial \psi}{\partial x} \tag{2}$$

Hence, $\phi(x, y)$ and $\psi(x, y)$ individually satisfy Laplace's equation

$$\frac{\partial^2 \phi}{\partial x^2} + \frac{\partial^2 \phi}{\partial y^2} = 0, \quad \frac{\partial^2 \psi}{\partial x^2} + \frac{\partial^2 \psi}{\partial y^2} = 0 \tag{3}$$

and it follows that any analytic function of a complex variable can be regarded as a complex potential distribution.

Numerous examples of complex potential distributions can be found in standard textbooks [1, 2, 3, 4]. Although they are mathematically rigorous, they only apply in idealised physical situations. In fluid mechanics, for example, they are restricted to the two-dimensional flow of an ideal fluid, which is assumed to be incompressible, inviscid and irrotational (i.e.

$\nabla \times \mathbf{u} = 0$, where $\mathbf{u}(x, y)$ is the velocity of the fluid). In this context, $\phi(x, y)$ and $\psi(x, y)$ are the 'velocity potential' and 'stream-function', respectively, and curves $\phi(x, y) = $ constant and $\psi(x, y) = $ constant are 'equipotentials' and 'stream-lines'.

From a mathematical perspective, however, complex potential distributions are not restricted to two dimensions. They can be analytically continued on a Riemann surface, as a set of parallel *x-y* planes, separated by branch cuts, or alternatively as an appropriate three-dimensional surface. In both representations, every point on the Riemann surface is intersected by an equipotential and a stream-line [5]. Interestingly, although Riemann surfaces are now regarded as abstract constructions for representing complex functions, they were originally conceived in terms of the flow of electric current between equipotentials and stream-lines in a thin conducting sheet [6].

This paper describes the potential flow of an ideal fluid on a three-dimensional Riemann surface corresponding to the Dirichlet eta function

$$\eta(s) = \sum_{n=1}^{\infty} \frac{(-1)^{n-1}}{n^s} \qquad \sigma > 0 \qquad (4)$$

[7,8], where $s = \sigma + it$ is a complex variable, and likewise the Dirichlet eta function

$$\eta(1-s) = \sum_{n=1}^{\infty} \frac{(-1)^{n-1}}{n^{1-s}} \qquad \sigma < 1 \qquad (5)$$

We focus on the behaviour of $\eta(s)$ and $\eta(1-s)$ in the 'critical strip', $0 < \sigma < 1$, and restrict attention to the semi-infinite interval, $0 \leq t < \infty$, since $\eta(s)$ and $\eta(1-s)$ are both symmetric about $t = 0$. $\eta(s)$ is closely related to the Riemann zeta function,

$$\zeta(s) = \sum_{n=1}^{\infty} \frac{1}{n^s} \qquad \sigma > 1 \qquad (6)$$

and it enables $\zeta(s)$ to be analytically continued in the critical strip, by the formula [7,8]

$$\zeta(s) = \frac{\eta(s)}{1 - 2^{1-s}} \qquad \sigma > 0 \qquad (7)$$

Riemann proved that $\zeta(s)$ has an infinite number of trivial zeros on the line $t = 0$, at $\sigma = -2, -4, -6, \ldots$, and an infinite number of non-trivial zeros in the critical strip [8, 9, 10, 11, 12]. He conjectured that all the zeros of $\zeta(s)$ in the critical strip are on the critical line, $\sigma = \frac{1}{2}$, which has since become known as the 'Riemann hypothesis'. The verification of the Riemann hypothesis would have significant consequences in pure mathematics [8]. In particular, it would imply a strong bound on the error of the prime number theorem by the inequality

$$|\pi(n) - \text{Li}(n)| < \frac{1}{8\pi}\sqrt{n}\ln n \qquad n \geq 2657 \qquad (8)$$

where $\pi(n)$ is the number of primes less than a given integer, $n$, and $\text{Li}(n)$ is the logarithmic integral [13, 14].

It follows from (7) that the zeros of $\eta(s)$ and $\zeta(s)$ are identical in the critical strip, so the Riemann hypothesis also applies to $\eta(s)$ [8]. It should also be noted from (7) that $\eta(s)$ has an infinite number of zeros on the line $\sigma = 1$, situated at

$$s_n = 1 \pm i\frac{2\pi n}{\ln 2} \qquad n = 1, 2, 3, \ldots \qquad (9)$$

but these zeros are irrelevant to the Riemann hypothesis. Furthermore, $\eta(s)$ is a holomorphic function whereas $\zeta(s)$ has a simple pole at $s = 1$, and this makes $\eta(s)$ somewhat easier to handle in the critical strip than $\zeta(s)$.

There is substantial body of circumstantial evidence in support of the Riemann hypothesis. In particular:

(i)  $\zeta(s)$ has an infinite number of zeros on the critical line [15];
(ii) at least two-fifths of the zeros of $\zeta(s)$ in the critical strip are on the critical line [16];
(iii) every zero of $\zeta(s)$ in the critical strip for $0 \leq t < 2.9 \times 10^9$ has a real part $\sigma = \frac{1}{2}$ [17].

This paper shows how potential theory provides information about the Dirichlet eta function that appears to have been overlooked in the number theory literature, in terms of properties that are mathematically rigorous and apply throughout the critical strip. These properties, together with well-known results from number theory, severely restrict the number of possible configurations of equipotentials and stream-lines near a zero in the critical strip, and we show how this leads to the verification of the Riemann hypothesis. [N.B. This paper replaces an earlier paper by the author [18], which did not use complex potential theory and did not consider all possible configurations of the Riemann surface.]

## 2    Riemann surfaces for $\eta(s)$ and $\eta(1-s)$

We define the $x$- and $y$-coordinates of a point $P(x_P, y_P, z_P)$ on the Riemann surface for

$$\eta(s) = \sum_{n=1}^{\infty}\frac{(-1)^{n-1}}{n^{\sigma+it}} = x_P(\sigma,t) + iy_P(\sigma,t) \qquad (10)$$

as

$$x_P(\sigma,t) = \sum_{n=1}^{\infty}\frac{(-1)^{n-1}}{n^{\sigma}}\cos(t\ln n), \quad y_P(\sigma,t) = -\sum_{n=1}^{\infty}\frac{(-1)^{n-1}}{n^{\sigma}}\sin(t\ln n) \qquad (11)$$

$x_P(\sigma,t)$ and $y_P(\sigma,t)$ are real functions of $\sigma$ and $t$, and the infinite series for $x_P(\sigma,t)$ and $y_P(\sigma,t)$ converge for $\sigma > 0$ [7].

Since $x_P(\sigma,t)$ and $y_P(\sigma,t)$ are multi-valued functions of $\sigma$ and $t$, we also assign a $z$-coordinate

$$z_P(t) = t \tag{12}$$

Hence, the Riemann surface for $\eta(s)$ in the critical strip is a three-dimensional surface which extends to infinity in the $z$-direction.

Likewise, we define the $x$- and $y$-coordinates of a point $Q(x_Q, y_Q, z_Q)$ on the Riemann surface for $\eta(1-s)$ by putting

$$\eta(1-s) = \sum_{n=1}^{\infty} \frac{(-1)^{n-1}}{n^{1-\sigma-it}} = x_Q(\sigma,t) + i y_Q(\sigma,t) \qquad \sigma < 1 \tag{13}$$

where

$$x_Q(\sigma,t) = \sum_{n=1}^{\infty} \frac{(-1)^{n-1}}{n^{1-\sigma}} \cos(t \ln n), \quad y_Q(\sigma,t) = \sum_{n=1}^{\infty} \frac{(-1)^{n-1}}{n^{1-\sigma}} \sin(t \ln n) \tag{14}$$

and

$$z_Q(t) = t \tag{15}$$

The shapes of the Riemann surfaces for $\eta(s)$ and $\eta(1-s)$ in the critical strip over the interval $0 \le t \le 26$ are shown in Figs. 1(a) and 1(b), respectively.

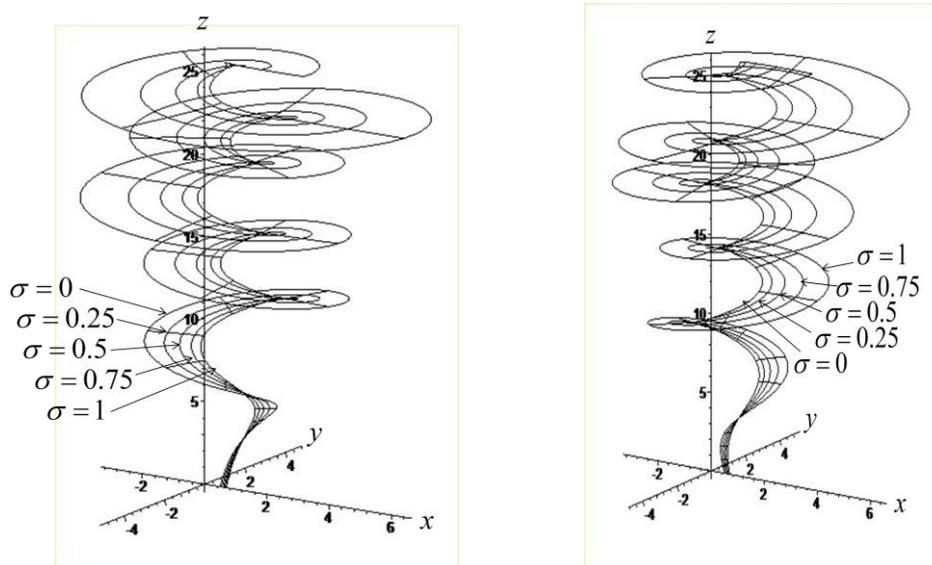

**Figs. 1(a),(b) Riemann surfaces in the critical strip for $\eta(s)$ and $\eta(1-s)$**

Note from Figs. 1(a) and 1(b) that:

(i)      curves $\sigma = \text{constant}$ on the Riemann surface for $\eta(s)$ are irregular spirals, and likewise for $\eta(1-s)$;

(ii)      any given curve $t = \text{constant}$ lies in the plane $z = t = \text{constant}$;

(iii)      the sequence of curves $\sigma = \text{constant}$ decreases from the innermost curve to the outermost curve on the Riemann surface for $\eta(s)$, but increases from the innermost curve to the outermost curve on the Riemann surface for $\eta(1-s)$;

(iv)      $\eta(s) = x_P(\sigma,t) + i y_P(\sigma,t) = 0$ whenever a curve $\sigma = \text{constant}$ intersects the $z$-axis.

The symmetry of the Riemann surfaces for $\eta(s)$ and $\eta(1-s)$ in the critical strip follows from the algebraic form of the individual terms in the series expansions for $\eta(s)$ and $\eta(1-s)$ in (11) and (14). Thus,

$$x_P(1-\sigma,t) = x_Q(\sigma,t), \quad y_P(1-\sigma,t) = -y_Q(\sigma,t) \qquad 0 < \sigma < 1 \qquad (16)$$

The spiral nature of the Riemann surface for $\eta(s)$ arises from the fact that particular curves $\sigma = \text{constant}$ are 'anchored' to the $z$-axis at the zeros of $\eta(s)$. From (9), the curve $\sigma = 1$ is anchored to the $z$-axis at the trivial zeros

$$t_n = \frac{2\pi n}{\ln 2} \qquad (n = 1, 2, 3, \ldots) \qquad (17)$$

and the curve $\sigma = \tfrac{1}{2}$ is anchored to the $z$-axis at an infinite number of non-trivial zeros $t \approx 14.1, 21.0, 25.0 \ldots$ [15]. Similarly, the spiral behaviour of $\eta(1-s)$ is due to the anchoring of the curves $\sigma = 0$ and $\sigma = \tfrac{1}{2}$ on the $z$-axis at the zeros of $\eta(1-s)$.

Figs. 2(a)-2(f) show the projections of curves $\sigma = \text{constant}$ and $t = \text{constant}$ on the Riemann surface for $\eta(s)$ onto the plane $z = 0$, for $0 \leq t \leq 11$, $11 \leq t \leq 15$, $15 \leq t \leq 18.5$, $18.5 \leq t \leq 21$, $21 \leq t \leq 24$, $24 \leq t \leq 26$. The curves $\sigma = \text{constant}$ and $t = \text{constant}$ are computed using MAPLE (and independently checked by evaluation of $x_P(\sigma,t)$ and $y_P(\sigma,t)$ from (11) using Aitken's $\Delta^2$ method). The minimum and maximum values of $t$ in any interval are chosen such that the curves $t = t_{min}$ and $t = t_{max}$ do not intersect one another. The corresponding projections for $\eta(1-s)$ are shown in Figs. 3(a)-3(f).

The small circles enclosing the origin in Figs. 2(b), 2(e) and 2(f) highlight the non-trivial zeros of $\eta(s)$, i.e. where the curve $\sigma = \tfrac{1}{2}$ intersects the origin at $t \approx 14.1, 21.0, 25.0$, respectively, and similarly for $\eta(1-s)$ in Figs. 3(b), 3(e) and 3(f). The small squares enclosing the origin in Figs. 2(a) and 2(c) highlight the trivial zeros of $\eta(s)$ on the line $\sigma = 1$, and likewise the trivial zeros of $\eta(1-s)$ on the line $\sigma = 0$.

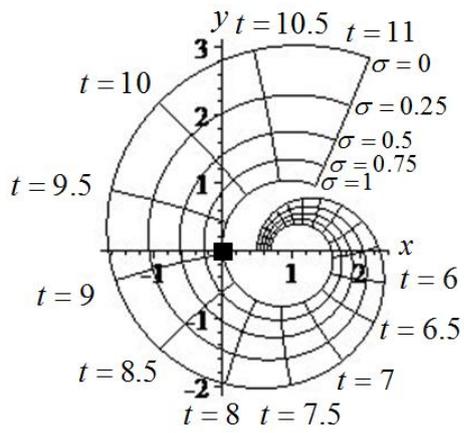

**Fig. 2a**

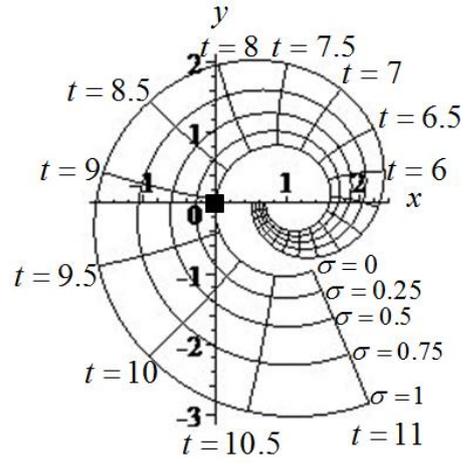

**Fig. 3a**

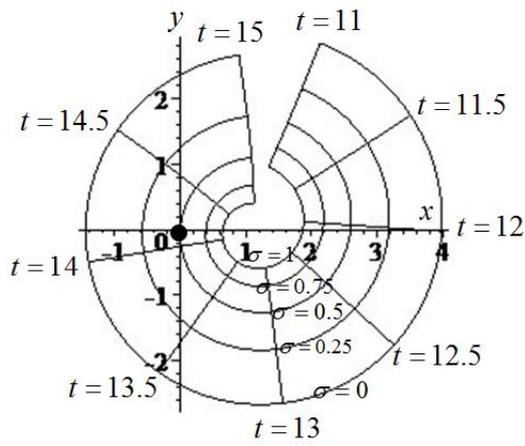

**Fig. 2b**

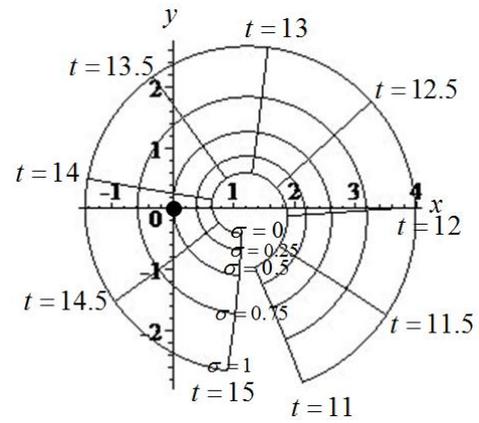

**Fig. 3b**

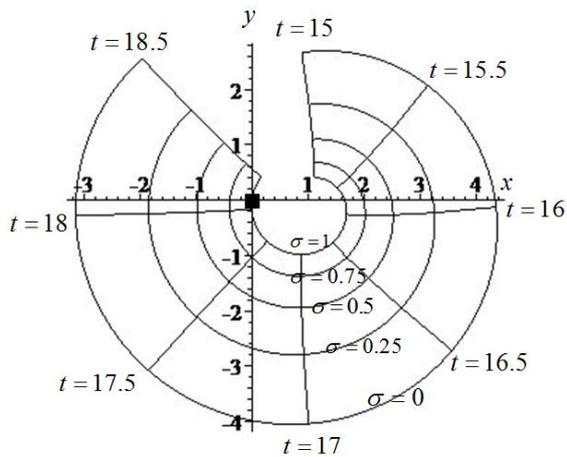

**Fig. 2c**

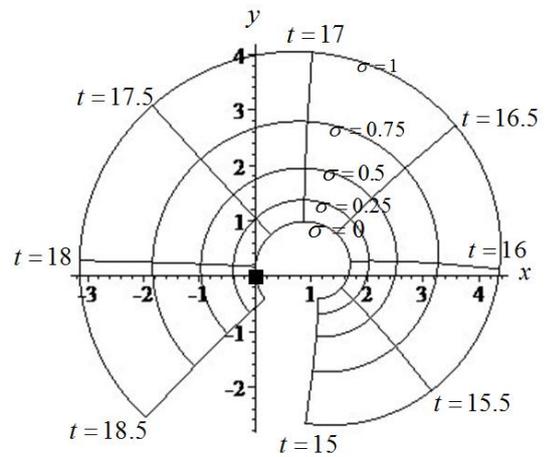

**Fig. 3c**

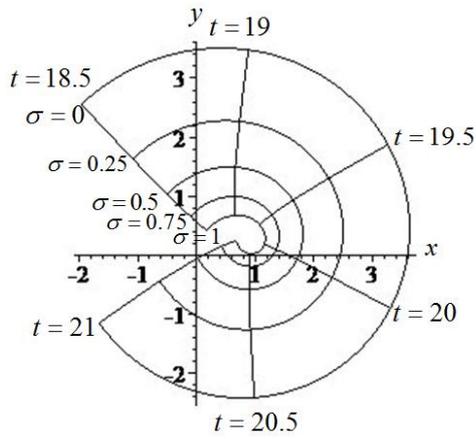

Fig. 2d

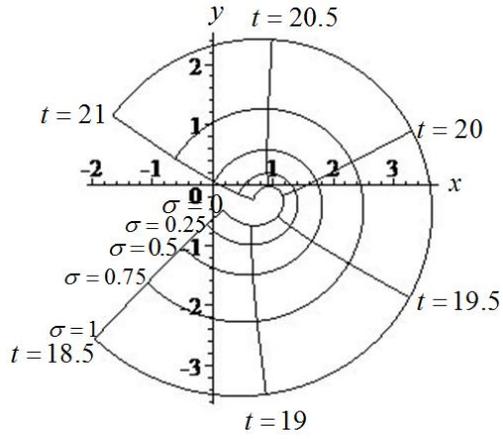

Fig. 3d

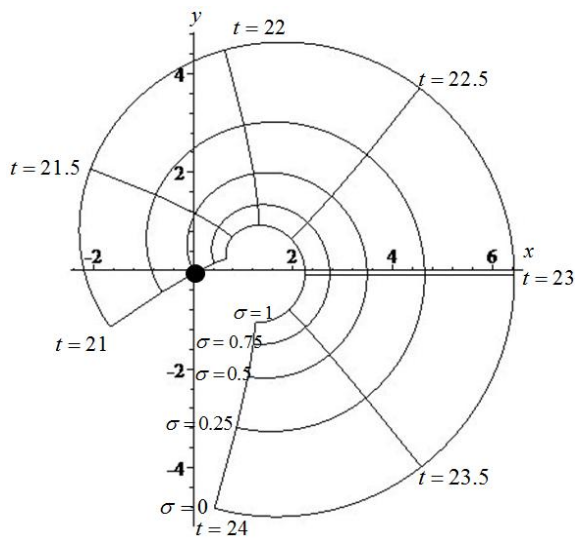

Fig. 2e

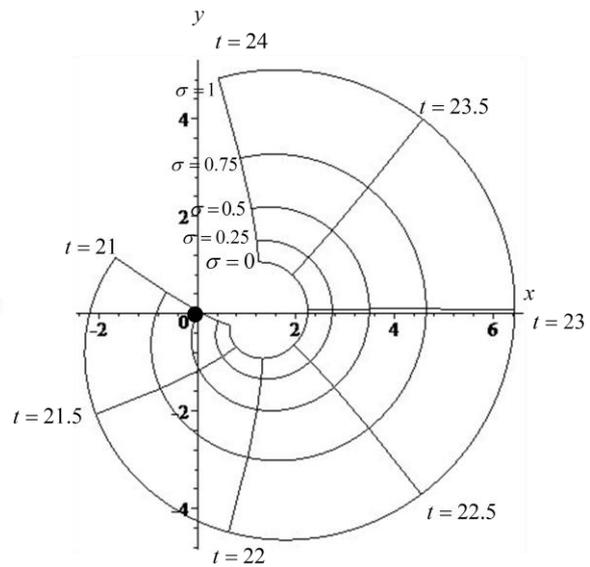

Fig. 3e

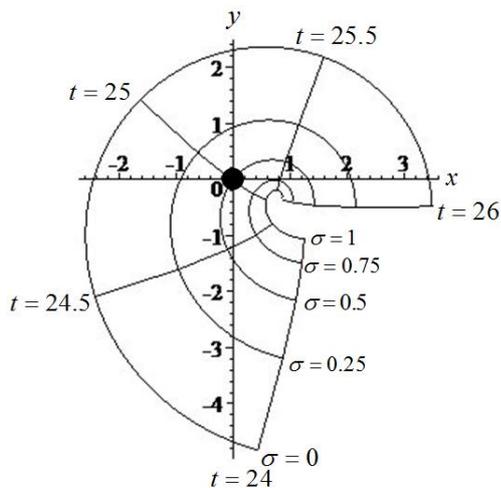

Fig. 2f

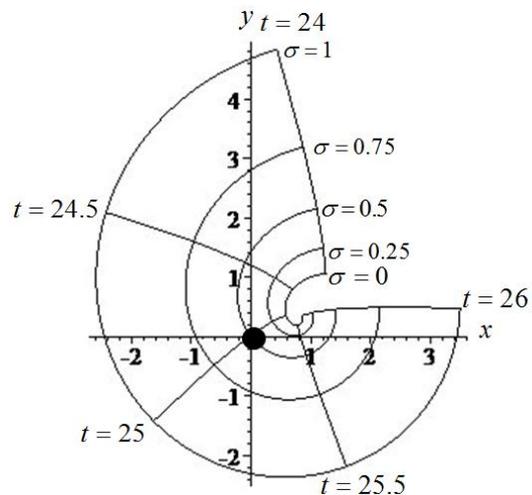

Fig. 3f

**Fig. 2a-2f**   Projection of curves $\sigma = $ const. and $t = $ const. corresponding to $\eta(s)$

**Fig. 3a-3f**   Projection of curves $\sigma = $ const. and $t = $ const. corresponding to $\eta(1-s)$

It is readily shown that any given point on the Riemann surface for $\eta(s)$, and all points in its neighbourhood, can be uniquely projected onto corresponding points in the $x$-$y$ plane. Thus, consider point $A(x,y,z)$ and its neighbourhood $D$, and the corresponding point $A'(x,y,0)$ and its neighbourhood $D'$ in the plane $z=0$, as shown in Fig. 4. $D$ is chosen to be sufficiently small that every curve $\sigma=$ constant is non-self-intersecting in $D$, and likewise every curve $t=$ constant (see Corollary to Theorem 1 in the next section).

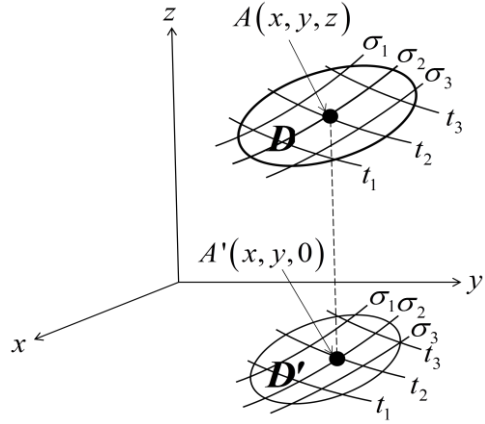

**Fig. 4  Projection of domain $D$ on the Riemann surface onto domain $D'$ in the $x$-$y$ plane**

Clearly, if the Riemann surface is not parallel to the $z$-axis at $A(x,y,z)$, or at any point in $D$, then all curves $\sigma=$ constant and $t=$ constant in $D$ can be uniquely projected onto $D'$.

Now suppose that a point $A_s(x,y,z)$ exists at which the Riemann surface is parallel to the $z$-axis, with Cartesian coordinates

$$x = x_P(\sigma,t), \quad y = y_P(\sigma,t), \quad z = t \tag{18}$$

Hence, $x_P(\sigma,t)$ and $y_P(\sigma,t)$ are stationary with respect to $t$ at $A_s(x,y,z)$, and

$$\frac{\partial x_P}{\partial t} = \frac{\partial y_P}{\partial t} = 0 \tag{19}$$

Also, from (11) and (19), we have

$$\frac{\partial x_P}{\partial \sigma} = \frac{\partial y_P}{\partial t} = 0 \tag{20}$$

and

$$\frac{\partial y_P}{\partial \sigma} = -\frac{\partial x_P}{\partial t} = 0 \tag{21}$$

at $A_s(x,y,z)$.

(20) and (21) imply that

$$\frac{d\eta}{ds} = \frac{\partial x_P}{\partial \sigma} + i\frac{\partial y_P}{\partial \sigma} = -i\frac{\partial x_P}{\partial t} + \frac{\partial y_P}{\partial t} = 0 \tag{22}$$

at $A_s(x, y, z)$. Since $\eta'(s)$ is a holomorphic function of $s$, and the zeros of a holomorphic function are *isolated*, it follows that $\eta'(s) \neq 0$ at all points in a (deleted) neighbourhood of $A_s(x, y, z)$. Hence, any point on the Riemann surface and all points in its neighbourhood can be uniquely projected onto a corresponding neighbourhood in the $x$-$y$ plane. This enables us to perform much of the following analysis in the $x$-$y$ plane.

## 3   Potential flow

We now derive a complex potential function, $\Omega(\eta)$, corresponding to $\eta(s)$. [N.B. We restrict attention to $\eta(s)$ in this section, since the same arguments apply to $\eta(1-s)$. Also, for simplicity, we omit the subscripts $P$ and $Q$, except where essential to the argument.]

The basic idea is to reverse the roles of $\eta$ and $s$, such that $s$ becomes a function of the complex variable $\eta = x + iy$, i.e.

$$s(\eta) = \sigma(x, y) + it(x, y) \tag{23}$$

where $\sigma(x, y)$ and $t(x, y)$ are real functions of $x$ and $y$. To justify this transformation, we prove

**Theorem 1**    $s(\eta)$ **is an analytic function of** $\eta$

Since $\eta(s)$ is a holomorphic function (i.e. it has no singularities), $\eta'(s)$ is also a holomorphic function. Since the zeros of a holomorphic function are isolated, the zeros of $\eta'(s)$ are isolated. Hence, the inverse derivative

$$\frac{ds}{d\eta} = \frac{\partial \sigma}{\partial x} + i\frac{\partial t}{\partial x} = \frac{\partial t}{\partial y} - i\frac{\partial \sigma}{\partial y} \tag{24}$$

exists, except (possibly) at isolated points where $\eta'(s) = 0$.

It follows from (24) that the partial derivatives of $\sigma(x, y)$ and $t(x, y)$ satisfy the Cauchy-Riemann conditions

$$\frac{\partial \sigma}{\partial x} = \frac{\partial t}{\partial y}, \qquad \frac{\partial \sigma}{\partial y} = -\frac{\partial t}{\partial x} \tag{25}$$

$s(\eta)$ is therefore an analytic function of $\eta$, except (possibly) at isolated points.    □

**Corollary: Curves $\sigma(x, y) = \text{constant}$ are non-self-intersecting in any sufficiently small element in the $x$-$y$ plane, and likewise curves $t(x, y) = \text{constant}$**

Consider a rectangular element of the $x$-$y$ plane with sides $\delta x$ and $\delta y$, where $(\delta x)^2 \ll 1$ and $(\delta y)^2 \ll 1$. Curves $\sigma(x, y) = \text{constant}$ in the element satisfy the equation

$$\delta\sigma = \frac{\partial \sigma}{\partial x}\delta x + \frac{\partial \sigma}{\partial y}\delta y + \varepsilon = 0 \tag{26}$$

where $\varepsilon$ is of order $(\delta x)^2, (\delta y)^2$. From Theorem 1, the coefficients of $\delta x$ and $\delta y$ in (26) are finite quantities. For sufficiently small values of $(\delta x)^2$ and $(\delta y)^2$, terms of order $\varepsilon$ can be ignored, and any given curve $\sigma(x, y) = \text{constant}$ tends to a straight line, of slope

$$\frac{dy}{dx} \approx -\frac{\partial \sigma/\partial x}{\partial \sigma/\partial y} \tag{27}$$

It follows that any curve $\sigma(x, y) = \text{constant}$ is non-self-intersecting inside the element. The same argument applies to any curve $t(x, y) = \text{constant}$. □

To be consistent with the conventions of fluid mechanics, we now define a complex potential function, $\Omega(\eta)$, such that curves $\sigma(x, y) = \text{constant}$ follow stream-lines, and curves $t(x, y) = \text{constant}$ follow equipotentials. Thus, we put

$$\Omega(\eta) = -is(\eta) = t(x, y) - i\sigma(x, y) \equiv \phi(x, y) + i\psi(x, y) \tag{28}$$

Hence

$$\phi(x, y) = t(x, y), \quad \psi(x, y) = -\sigma(x, y) \tag{29}$$

It follows from (25) and (29) that $\phi(x, y)$ and $\psi(x, y)$ satisfy the Cauchy-Riemann conditions

$$\frac{\partial \phi}{\partial x} = \frac{\partial \psi}{\partial y}, \quad \frac{\partial \phi}{\partial y} = -\frac{\partial \psi}{\partial x} \tag{30}$$

Eliminating $\psi(x, y)$ yields Laplace's equation

$$\frac{\partial^2 \phi}{\partial x^2} + \frac{\partial^2 \phi}{\partial y^2} = 0 \tag{31}$$

Similarly, eliminating $\phi(x, y)$ we obtain

$$\frac{\partial^2 \psi}{\partial x^2} + \frac{\partial^2 \psi}{\partial y^2} = 0 \tag{32}$$

Also note, from (12) and (29), that

$$z = t(x, y) = \phi(x, y) \tag{33}$$

Hence, the equipotential $\phi(x, y) = $ constant follows the curve $t(x, y) = $ constant at a height $z$ above the $x$-$y$ plane.

Differentiating $\Omega(\eta)$ with respect to $\eta$, using (28) and (30), we obtain

$$\frac{d\Omega}{d\eta} = -i\frac{ds}{d\eta} = \frac{\partial \phi}{\partial x} + i\frac{\partial \psi}{\partial x} = \frac{\partial \psi}{\partial y} - i\frac{\partial \phi}{\partial y} \equiv u(x, y) - iv(x, y) \tag{34}$$

where

$$u(x, y) = \frac{\partial \phi}{\partial x} = \frac{\partial \psi}{\partial y} \tag{35}$$

and

$$v(x, y) = \frac{\partial \phi}{\partial y} = -\frac{\partial \psi}{\partial x} \tag{36}$$

are the components of velocity of the fluid in the $x$- and $y$-directions, respectively [1].

Also, from (33) we have $\delta z = \delta\phi(x, y)$, so the component of velocity in the $z$-direction is simply given by

$$w = \frac{d\phi}{dz} \equiv 1, \text{ independent of } x \text{ and } y. \tag{37}$$

We now derive three theorems relating to the potential flow field, which are used in the following section to show that zeros of $\eta(s)$ and $\eta(1-s)$ cannot exist in the critical strip off the critical line. [N.B. These theorems are well-known in potential theory but do not seem to have been used in relation to the Riemann hypothesis, and are included here for completeness. Also note that the following proofs do not involve any assumptions of an empirical nature, in order to maintain mathematical rigour.]

**Theorem 2    Stream-lines and equipotentials are mutually orthogonal**

Expanding $\psi(x, y)$ to first order in $\delta x$ and $\delta y$ along a stream-line $\psi(x, y) = $ constant (see Fig. 5), and using (35) and (36), we have

$$\delta\psi(x,y) = \frac{\partial \psi}{\partial x}\delta x + \frac{\partial \psi}{\partial y}\delta y = -v(x,y)\delta x + u(x,y)\delta y = 0 \qquad (38)$$

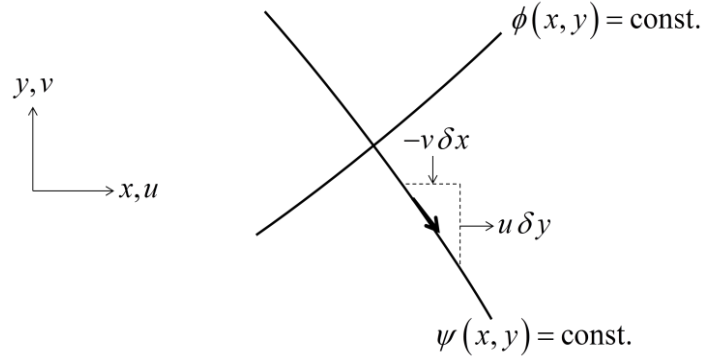

**Fig. 5**      **Flow parallel to the stream-line $\psi(x,y) =$ constant and orthogonal to the equipotential $\phi(x,y) =$ constant**

Hence, the gradient of the streamline $\psi(x,y) =$ constant is given by

$$\left(\frac{dy}{dx}\right)_{\psi=\text{const.}} = \frac{v(x,y)}{u(x,y)} \qquad (39)$$

It follows from (39) that the total velocity is parallel to the stream-line $\psi(x,y) =$ constant. Also, along the equipotential $\phi(x,y) =$ constant, we have

$$\delta\phi(x,y) = \frac{\partial \phi}{\partial x}\delta x + \frac{\partial \phi}{\partial y}\delta y = u(x,y)\delta x + v(x,y)\delta y = 0 \qquad (40)$$

The gradient of the equipotential $\phi(x,y) =$ constant is therefore given by

$$\left(\frac{dy}{dx}\right)_{\phi=\text{const.}} = -\frac{u(x,y)}{v(x,y)} \qquad (41)$$

(39) and (41) imply that

$$\left(\frac{dy}{dx}\right)_{\psi=\text{const.}}\left(\frac{dy}{dx}\right)_{\phi=\text{const.}} = -1 \qquad (42)$$

at the point of intersection of the curves $\psi(x,y) =$ constant and $\phi(x,y) =$ constant. Hence, the curves intersect one another orthogonally, and the transformation from the $s$-plane to the $x$-$y$ plane is conformal.   □

**Theorem 3**   The integral $\int \left[ -v(x,y)dx + u(x,y)dy \right]$ between any two stream-lines is invariant throughout the critical strip

Consider the closed region $S$ of the Riemann surface for $\eta(s)$ bounded by the contour $A_0 B_0 BAA_0$ shown in Fig. 6.

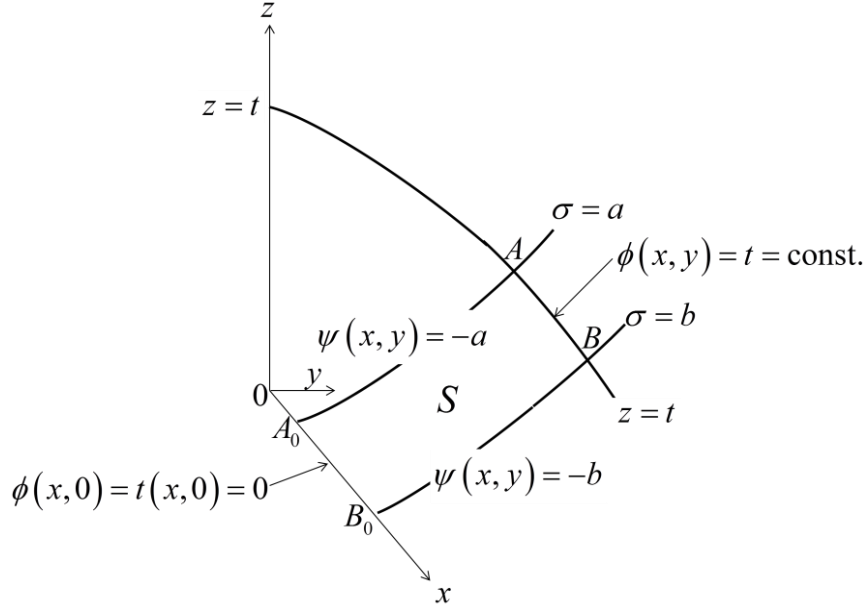

**Fig. 6  Flow between equipotentials and stream-lines**

Along the sides of the contour $A_0 B_0 BAA_0$ we put

$\phi(x,0) = t(x,0) = 0$   on $A_0 B_0$ (equipotential)
$\psi(x,y) = -b$   on $B_0 B$ (stream-line)
$\phi(x,y) = t = $ constant   on $BA$ (equipotential)
$\psi(x,y) = -a$   on $AA_0$ (stream-line)

where $0 < a < b < 1$.

By Green's theorem, if $F(x,y)$ and $G(x,y)$ and their first order partial derivatives are continuous functions of $x$ and $y$ in $S$ and on the contour $A_0 B_0 BAA_0$, then

$$\oint_{A_0 B_0 BAA_0} (Fdx + Gdy) = \iint_S \left( \frac{\partial G}{\partial x} - \frac{\partial F}{\partial y} \right) dx\, dy \tag{43}$$

Using (35) and (36), we put

$$F(x,y) = -v(x,y) = -\frac{\partial \phi}{\partial y} = \frac{\partial \psi}{\partial x}, \qquad G(x,y) = u(x,y) = \frac{\partial \phi}{\partial x} = \frac{\partial \psi}{\partial y} \tag{44}$$

Substituting for $F(x, y)$ and $G(x, y)$ in (43), and using Laplace's equation (31), we obtain

$$\oint_{A_0B_0BAA_0}[-v(x,y)dx+u(x,y)dy] = \iint_S\left(\frac{\partial^2\phi}{\partial x^2}+\frac{\partial^2\phi}{\partial y^2}\right)dx\,dy = \iint_S\left(\frac{\partial^2\psi}{\partial x\partial y}-\frac{\partial^2\psi}{\partial y\partial x}\right)dx\,dy = 0$$

Separating the integrals along each segment,

$$\left(\int_{A_0}^{B_0}+\int_{B_0}^{B}+\int_{B}^{A}+\int_{A}^{A_0}\right)[-v(x,y)dx+u(x,y)dy] = 0 \tag{45}$$

From (38), the integrals along the stream-lines $B_0B$ and $A_0A$ in (45) are identically zero. Also, along the equipotential $A_0B_0$, we have $\phi(x,0) = t(x,0) \equiv 0$ $(a \leq x \leq b)$, so that $u(x,0) = \partial\phi/\partial x \equiv 0$.

(45) reduces to

$$\int_{A_0}^{B_0}-v(x,0)dx + \int_{B}^{A}[-v(x,y)dx+u(x,y)dy] = 0 \tag{46}$$

From (36), the first integral in (46) is given by

$$\int_{A_0}^{B_0}-v(x,0)dx = \int_{A_0}^{B_0}\left(\frac{\partial\psi}{\partial x}\right)_{x,0}dx = \psi_{B_0} - \psi_{A_0} = -b-(-a) = -b+a \tag{47}$$

Substituting in (46), we obtain

$$\int_{B}^{A}[-v(x,y)dx+u(x,y)dy] = b-a \tag{48}$$

Since $b$ and $a$ are constants, the integral on the left hand side of (48) is invariant throughout the critical strip, i.e. independent of $t$, and the theorem is proved.  □

[N.B. In fluid mechanics, the integral in (48) represents the mass flux between any two stream-lines. The theorem is usually derived in textbooks of fluid mechanics by assuming mass conservation in every element of the fluid. However, the use of Green's theorem (43) in the above derivation avoids the need to make this assumption.]

**Theorem 4**    $\sigma(x, y)$ **varies monotonically along any curve** $t(x, y) = \text{constant}$

Suppose that the theorem is false, by assuming that $\sigma(x, y)$ is a non-monotonic function along some (equipotential) curve $\phi(x, y) = t(x, y) = \text{constant}$. It follows that there is at least one stationary point, $A_s(x, y)$, on the curve, where

$$\delta\sigma(x, y) = \frac{\partial \sigma}{\partial x}\delta x + \frac{\partial \sigma}{\partial y}\delta y = 0 \qquad (49)$$

Putting $\sigma(x, y) = -\psi(x, y)$ from (29), and using (35) and (36), (49) becomes

$$-\delta\psi(x, y) = -\frac{\partial \psi}{\partial x}\delta x - \frac{\partial \psi}{\partial y}\delta y = v(x, y)\delta x - u(x, y)\delta y = 0 \qquad (50)$$

Also, along the equipotential $\phi(x, y) = t(x, y) = \text{constant}$, we have

$$\delta\phi(x, y) = \frac{\partial \phi}{\partial x}\delta x + \frac{\partial \phi}{\partial y}\delta y = u(x, y)\delta x + v(x, y)\delta y = 0 \qquad (51)$$

Eliminating $\delta x$ and $\delta y$ between (50) and (51) yields

$$\frac{u(x, y)}{v(x, y)} = -\frac{v(x, y)}{u(x, y)}$$

or

$$[u(x, y)]^2 = [v(x, y)]^2 = 0 \qquad (52)$$

Hence

$$u(x, y) = v(x, y) = 0 \qquad (53)$$

at $A_s(x, y)$.

(34) and (53) imply that

$$-i\frac{ds}{d\eta} = u(x, y) - iv(x, y) = 0 \qquad (54)$$

at $A_s(x, y)$. Inverting (54), we have

$$\frac{d\eta}{ds} = \frac{i}{u(x, y) - iv(x, y)} = \frac{i}{0} \qquad (55)$$

at $A_s(x, y)$. Hence, $\eta'(s)$ is non-analytic at $A_s(x, y)$, implying that $\eta'(s)$ is not a holomorphic function of $s$. The assumption that $\sigma(x, y)$ is non-monotonic at some point along a curve $t(x, y) = \text{constant}$ is therefore false and the theorem is proved. □

## 4    Zeros in the critical strip

We now focus on the central issue of whether zeros of $\eta(s)$ can exist in the critical strip off the critical line, $\sigma = \frac{1}{2}$. We begin by proving

**Theorem 5**    If $\eta(s) = 0$ at some point $s = \sigma^* + it^*$ $\left(\sigma^* \neq \frac{1}{2}\right)$ in the critical strip, then $\eta(s) = 0$ at the point $s = (1-\sigma^*) + it^*$, also in the critical strip

Suppose that $\eta(s) = 0$ at some point $s^* = \sigma^* + it^*$, where $0 < \sigma^* < \frac{1}{2}$ and $0 \leq t^* < \infty$. It follows from (10) that

$$\eta(s^*) = x_P(\sigma^*, t^*) + i y_P(\sigma^*, t^*) = 0 \tag{56}$$

Hence

$$x_P(\sigma^*, t^*) = y_P(\sigma^*, t^*) = 0 \tag{57}$$

Also, $\eta(s)$ and $\eta(1-s)$ satisfy the functional equation

$$\eta(1-s) = k(s)\,\eta(s) \tag{58}$$

[19], where

$$k(s) = -\pi^{-s}\left(\frac{2^s - 1}{2^{s-1} - 1}\right)\cos\left(\tfrac{1}{2}\pi s\right)\Gamma(s) \tag{59}$$

Taking the modulus of both sides of (59), we have

$$|k(s)| = |\pi^{-s}|\left|\frac{2^s - 1}{2^{s-1} - 1}\right|\left|\cos\left(\tfrac{1}{2}\pi s\right)\right||\Gamma(s)| \tag{60}$$

Since

(i)    $\left|\pi^{-s}\right| = \pi^{-\sigma} \neq 0$

(ii)    $\left|\dfrac{2^s - 1}{2^{s-1} - 1}\right| \neq 0$, for $0 < \sigma < 1$

(iii)    $\left|\cos\left(\tfrac{1}{2}\pi s\right)\right| = \sqrt{\cos^2\left(\tfrac{1}{2}\pi\sigma\right)\cosh^2\left(\tfrac{1}{2}\pi t\right) + \sin^2\left(\tfrac{1}{2}\pi\sigma\right)\sinh^2\left(\tfrac{1}{2}\pi t\right)} \neq 0$

(iv)    $\Gamma(s) \neq 0$ [20]

it follows that $|k(s)| \neq 0$ in the interval $0 < \sigma < 1$.

(13) and (58) then imply that

$$\eta(1-s^*) = x_Q(\sigma^*, t^*) + i\, y_Q(\sigma^*, t^*) = 0 \tag{61}$$

Hence

$$x_Q(\sigma^*, t^*) = y_Q(\sigma^*, t^*) = 0 \tag{62}$$

Also, from the symmetry relations (16), we have

$$x_P(1-\sigma^*, t^*) = x_Q(\sigma^*, t^*), \quad y_P(1-\sigma^*, t^*) = -y_Q(\sigma^*, t^*) \tag{63}$$

and it follows from (62) and (63) that

$$x_P(1-\sigma^*, t^*) = y_P(1-\sigma^*, t^*) = 0$$

Hence,

$$\eta_P(1-\sigma^* + it^*) = x_P(1-\sigma^*, t^*) + i\, y_P(1-\sigma^*, t^*) = 0 \tag{64}$$

i.e. $\eta(s) = 0$ at $s = 1 - \sigma^* + it^*$, and the theorem is proved. □

Theorem 2 and Theorem 5 impose a severe restriction on the number of possible configurations of the curves $\sigma(x, y) = \sigma^*$, $\sigma(x, y) = 1 - \sigma^*$ ($0 < \sigma^* < \tfrac{1}{2}$), and $t(x, y) = t^*$, in the neighbourhood of a zero of $\eta(s)$. The only situations in which both theorems are satisfied are illustrated (schematically) in Figs. 7(a) and 7(b) below.

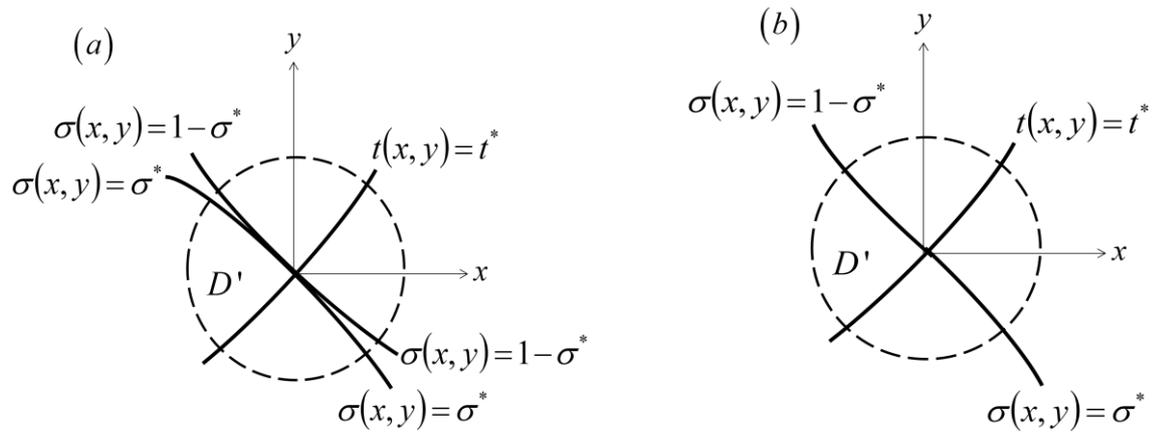

**Fig. 7(a),(b)   Configurations of the curves $\sigma(x, y) = \sigma^*, \sigma(x, y) = 1 - \sigma^*, t(x, y) = t^*$ satisfying Theorem 2 and Theorem 5**

In Fig. 7(a), curves $\sigma(x, y) = \sigma^*$ and $\sigma(x, y) = 1 - \sigma^*$ are both orthogonal to the curve $t(x, y) = t^*$ (Theorem 2) and meet at a tangent at the zero of $\eta(s)$ at $(x, y) = (0,0)$.

In Fig. 7(b), curves $\sigma(x, y) = \sigma^*$ and $\sigma(x, y) = 1 - \sigma^*$ lie on opposite sides of the curve $t(x, y) = t^*$ at $(x, y) = (0,0)$. Hence, if $\sigma^* \neq \tfrac{1}{2}$, then either  (i) $\sigma(x, y)$ is discontinuous at $(x, y) = (0, 0)$, or  (ii) the curve $t(x, y) = t^*$ follows a branch cut, which intersects $(x, y) = (0, 0)$.

(i) is rejected since $\sigma(x, y)$ is an analytic function of $x$ and $y$. However, (ii) cannot be immediately ruled out, because it is conceivable that $\sigma(x, y)$ could vary continuously across a branch cut from an upper $x$-$y$ plane to a lower $x$-$y$ plane, as shown in Figs. 8(a) and 8(b).

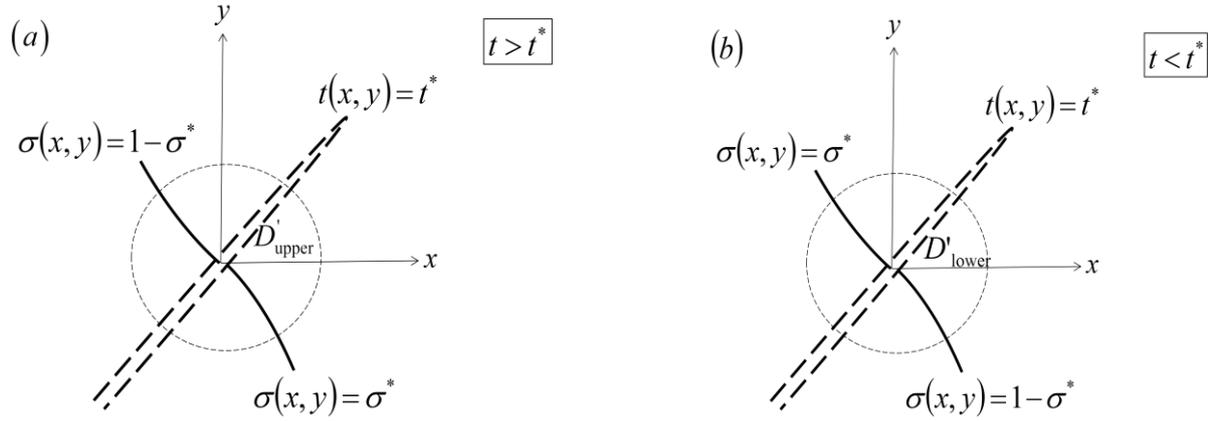

**Fig. 8(a), 8(b)  Curves $\sigma(x, y) = \sigma^*$ and $\sigma(x, y) = 1 - \sigma^*$ in two $x$-$y$ planes**

We now analyse all possible configurations of curves $\sigma(x, y) = \sigma^*$, $\sigma(x, y) = 1 - \sigma^*$, and $t(x, y) = t^*$ in the neighbourhood of a zero of $\eta(s)$, and show in each case that inconsistencies arise that contradict the assumption that zeros of $\eta(s)$ can exist in the critical strip for $\sigma^* \neq \tfrac{1}{2}$.
.

**I**  **Curves $\sigma(x, y) = \sigma^*$ and $\sigma(x, y) = 1 - \sigma^*$ meet at a tangent at a zero of $\eta(s)$**

We first show that the configuration shown in Fig. 7(a) cannot arise, by exploiting the invariance of the 'flux' integral in (48).

Consider the curves $\sigma(x, y) = \sigma^*$, $\sigma(x, y) = 1 - \sigma^*$ ($0 < \sigma^* < \tfrac{1}{2}$) and $t(x, y) = t^*$ ($0 \leq t^* < \infty$) shown in Figs. 9(a) and 9(b). Curves $\sigma(x, y) = \sigma^*$ and $\sigma(x, y) = 1 - \sigma^*$ intersect the curve $t(x,0) = 0$ (which follows the $x$-axis) at points $A_0$ and $B_0$, respectively, and intersect the curve $t(x, y) = t^*$ at points $A(\tfrac{1}{2}\delta x, \tfrac{1}{2}\delta y)$ and $B(-\tfrac{1}{2}\delta x, -\tfrac{1}{2}\delta y)$, respectively, on opposite sides of some point $Z(0,0, t^*)$ on the $z$-axis where $\eta(s) = x + iy = 0$.

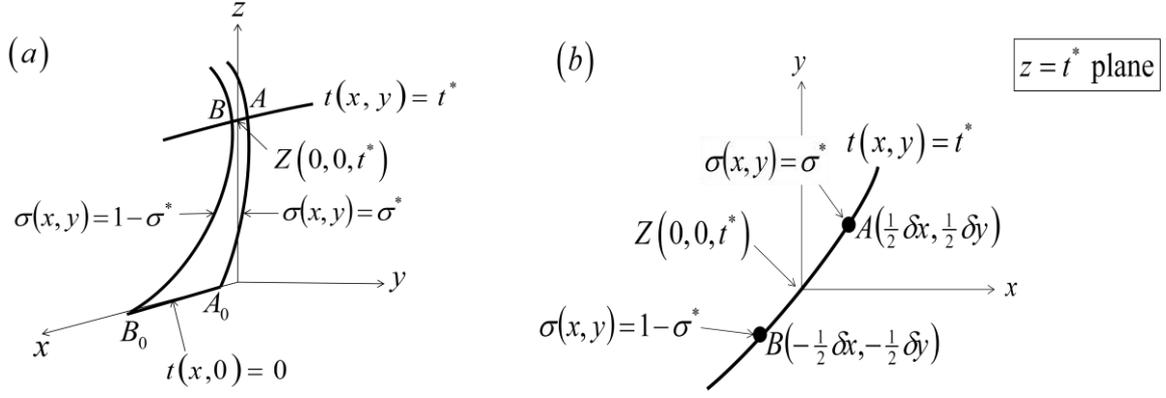

**Fig. 9(a), 9(b) Intersection of curves $\sigma(x, y) = \sigma^*$ and $\sigma(x, y) = 1 - \sigma^*$ with curve $t = t^*$**

From (48), we have

$$\int_B^A \left[ -v(x, y)dx + u(x, y)dy \right] = (1 - \sigma^*) - \sigma^* = 1 - 2\sigma^* \qquad (65)$$

For $(\delta x)^2 \ll 1$ and $(\delta y)^2 \ll 1$, (65) can be expressed in the form

$$-v(0,0)\delta x + u(0,0)\delta y + \varepsilon = 1 - 2\sigma^* \qquad (66)$$

where $\varepsilon$ is of order $(\delta x)^2, (\delta y)^2$.

For $0 < \sigma^* < \tfrac{1}{2}$, there are no solutions of (66) in the limit as $\delta x \to 0$ and $\delta y \to 0$ unless $u(0,0)$ and/or $v(0,0)$ are singular. Also, from (34), the expression

$$-i\frac{ds}{d\eta} = u(0,0) - iv(0,0) \qquad (67)$$

is singular, at $\eta = x + iy = 0$.

Using Laurent's theorem, we now expand $ds/d\eta$ in an annular neighbourhood of $\eta = 0$ in terms of the series

$$\frac{ds}{d\eta} = \sum_{n=0}^{\infty} a_n \eta^n + \frac{b_1}{\eta} + \frac{b_2}{\eta^2} + \ldots \qquad (68)$$

where $a_n$ and $b_n$ are complex constants. Integrating both sides of (68) with respect to $\eta$, we obtain

$$s(\eta) = \sum_{n=0}^{\infty} \frac{a_n}{n+1} \eta^{n+1} + b_1 \ln \eta - \frac{b_2}{\eta} - \ldots + C \qquad (69)$$

where $C$ is a constant of integration.

The right hand side of (69) is singular at $\eta = 0$, which implies that $s(0) = \sigma^* + it^*$ is not a finite quantity. However, by definition, $\sigma^*$ and $t^*$ lie in the intervals $0 < \sigma^* < \frac{1}{2}$ and $0 < t^* < \infty$, so the assumption that a zero of $\eta(s)$ exists for $\sigma^* \neq \frac{1}{2}$ must be false.

$\sigma^* = \frac{1}{2}$ is a special case, in which (66) reduces to the equation

$$-v(0,0)\delta x + u(0,0)\delta y + \varepsilon = 0 \qquad (70)$$

(70) can be satisfied for finite values of $u(0,0)$ and $v(0,0)$ in the limit as $\delta x \to 0$ and $\delta y \to 0$, so there is no contradiction in this case.

## II    Curves $\sigma(x, y) = \sigma^*$ and $\sigma(x, y) = 1 - \sigma^*$ intersect a branch cut

We next consider the possibility that the curves $\sigma(x, y) = \sigma^*$ and $\sigma(x, y) = 1 - \sigma^*$ ($\sigma^* \neq \frac{1}{2}$) intersect a branch cut which follows the curve $t(x, y) = t^*$, as shown in Figs. 8(a) and 8(b). In order for a branch cut to arise, the Riemann surface would need to intersect itself along part, or all, of the curve $z = t(x, y) = t^*$. We now analyse the following types of branch cut:

(i)    open-ended,
(ii)   half-open-ended,
(iii)  closed,

since all branch cuts can be classified according to one of these three types [21]. In each case we show that it is not possible to define a self-consistent branch cut.

### Case II(i)    Open-ended branch cut

First consider an open-ended branch cut which follows the curve $t(x, y) = t^*$, and is intersected by all curves $\sigma(x, y) = $ constant ($-\infty < \sigma < \infty$). Since $\sigma(x, y)$ is monotonic along $t(x, y) = t^*$ (Theorem 4), the extreme curves, $\sigma(x, y) = \infty$ and $\sigma(x, y) = -\infty$ intersect the open ends of the branch cut, at $E_\infty$ and $E_{-\infty}$, respectively, as shown in Fig. 10.

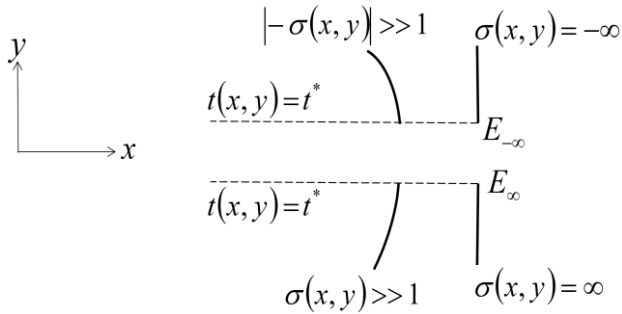

**Fig. 10**    Curves $\sigma(x, y) = $ const. near the open end of the branch cut (schematic)

From (11), the real part of $\eta(s)$ has upper and lower bounds given by

$$x_{upper}(\sigma, t^*) = 1 + \frac{1}{2^\sigma} + \frac{1}{3^\sigma} + \ldots = \zeta(\sigma) \tag{71}$$

$$x_{lower}(\sigma, t^*) = 1 - \frac{1}{2^\sigma} - \frac{1}{3^\sigma} - \ldots = 2 - \left(1 + \frac{1}{2^\sigma} + \frac{1}{3^\sigma} + \ldots\right) = 2 - \zeta(\sigma) \tag{72}$$

where

$$\zeta(\sigma) = \sum_{n=1}^{\infty} \frac{1}{n^\sigma} \qquad \sigma > 1 \tag{73}$$

is the Euler zeta function (i.e. the Riemann zeta function (6) on the line $t = 0$). It follows from (73) that

$$\zeta(\sigma) \to 1 \quad \text{as} \quad \sigma \to \infty \tag{74}$$

since every term in the series (73) tends to zero except the first term. Also, for $\sigma > 1$, $\zeta(\sigma)$ is a monotonically decreasing function of $\sigma$, as shown in Fig. 11.

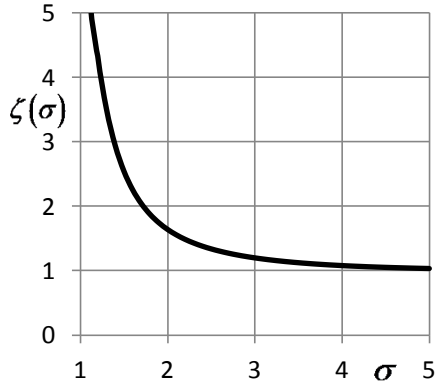

**Fig. 11**　　Euler zeta function

Similarly, from (11), the imaginary part of $\eta(s)$ has upper and lower bounds

$$y_{upper}(\sigma, t^*) = \frac{1}{2^\sigma} + \frac{1}{3^\sigma} + \ldots = \zeta(\sigma) - 1 \tag{75}$$

$$y_{lower}(\sigma, t^*) = -\frac{1}{2^\sigma} - \frac{1}{3^\sigma} - \ldots = 1 - \zeta(\sigma) \tag{76}$$

(71)-(72) and (75)-(76) imply that the coordinates of any point intersected by a curve $\sigma(x, y) \gg 1$ on the lower side of the branch cut in Fig. 10 are bounded by the inequalities

$$2 - \zeta(\sigma) < x(\sigma, t^*) < \zeta(\sigma), \quad 1 - \zeta(\sigma) < y(\sigma, t^*) < \zeta(\sigma) - 1 \tag{77}$$

Hence, for $\sigma \gg 1$,

$$x(\sigma, t^*) \approx 1, \quad y(\sigma, t^*) \approx 0 \tag{78}$$

Now consider the point of intersection of a curve $|-\sigma(x, y)| \gg 1$ with the curve $t(x, y) = t^*$ on the upper side of the branch cut in Fig. 10. Given that $\eta(s)$ has no zeros in the critical strip off the critical line for $0 < t^* < 2.9 \times 10^9$ [17], we can restrict attention to $t^* \gg 1$, which allows us to exploit well-known asymptotic results.

We first derive an asymptotic approximation for $|k(\sigma + i t^*)|$ for $\sigma \gg 1$, $t^* \gg 1$. Each term in (60) is given by

(i) $$\left| \frac{1}{\pi^{\sigma + i t^*}} \right| = \frac{1}{\pi^\sigma} \tag{79}$$

(ii) $$\left| \frac{2^{\sigma + i t^*} - 1}{2^{\sigma + i t^* - 1} - 1} \right| = \left| \frac{2^{i t^*} - 2^{-\sigma}}{2^{-1} 2^{i t^*} - 2^{-\sigma}} \right| \approx 2 \tag{80}$$

(iii) $$\left| \cos\left(\tfrac{1}{2} \pi (\sigma + i t^*)\right) \right| = \tfrac{1}{2} \left| e^{i \tfrac{1}{2} \pi (\sigma + i t^*)} + e^{-i \tfrac{1}{2} \pi (\sigma + i t^*)} \right| \approx \tfrac{1}{2} e^{\tfrac{1}{2} \pi t^*} \tag{81}$$

(iv) $$|\Gamma(\sigma + i t^*)| \geq \frac{\Gamma(\sigma)}{\sqrt{\cosh(\pi t^*)}} \approx 2 \pi^{\tfrac{1}{2}} e^{-\tfrac{1}{2} \pi t^*} e^{-\sigma} \sigma^{\sigma - \tfrac{1}{2}} \quad [20] \tag{82}$$

Substituting from (i)-(iv) above in (60), it follows that

$$|k(\sigma + i t^*)| \geq 2 e^{-\tfrac{1}{2}} \left( \frac{\sigma}{\pi e} \right)^{\sigma - \tfrac{1}{2}} \qquad \sigma \gg 1,\ t^* \gg 1 \tag{83}$$

Also, from (78), we have

$$|\eta(\sigma + i t^*)| = \sqrt{x^2(\sigma, t^*) + y^2(\sigma, t^*)} \approx 1 \qquad \sigma \gg 1 \tag{84}$$

Taking the modulus of both sides of (58), and using (83) and (84), then

$$|\eta(1 - \sigma - i t^*)| \geq 2 e^{-\tfrac{1}{2}} \left( \frac{\sigma}{\pi e} \right)^{\sigma - \tfrac{1}{2}} \qquad \sigma \gg 1,\ t^* \gg 1 \tag{85}$$

Moreover, since $|\eta(1 - \sigma - i t^*)| = |\eta(1 - \sigma + i t^*)|$ and $1 - \sigma \approx -\sigma$ for $\sigma \gg 1$,

$$|\eta(-\sigma + i t^*)| \geq 2 e^{-\tfrac{1}{2}} \left( \frac{\sigma}{\pi e} \right)^{\sigma - \tfrac{1}{2}} \qquad \sigma \gg 1,\ t^* \gg 1 \tag{86}$$

Hence

$$\sqrt{x^2(-\sigma, t^*) + y^2(-\sigma, t^*)} \geq 2e^{-\frac{1}{2}}\left(\frac{\sigma}{\pi e}\right)^{\sigma-\frac{1}{2}} \gg 1 \qquad \sigma \gg 1,\ t^* \gg 1 \qquad (87)$$

(84) and (87) imply that

$$\sqrt{x^2(-\sigma, t^*) + y^2(-\sigma, t^*)} \geq 2e^{-\frac{1}{2}}\left(\frac{\sigma}{\pi e}\right)^{\sigma-\frac{1}{2}} \gg \sqrt{x^2(\sigma, t^*) + y^2(\sigma, t^*)} \approx 1 \qquad (88)$$

for $\sigma \gg 1,\ t^* \gg 1$.

In the limit as $\sigma \to \infty$, $\sqrt{x^2(-\sigma, t^*) + y^2(-\sigma, t^*)} \to \infty$ and $\sqrt{x^2(\sigma, t^*) + y^2(\sigma, t^*)} \to 1$. It follows that points $E_\infty$ and $E_{-\infty}$ cannot lie on opposite sides at the open end of the branch cut, and the possibility of an open-ended branch cut is therefore rejected.

**Case II(ii)   Half-open-ended branch cut**

We next consider a half-open-ended branch cut, as shown in Fig. 12, whereby only the curves $\sigma(x, y) = \text{constant}$ in the interval $\sigma_0 \leq \sigma < \infty$ intersect the branch cut, where $\sigma_0$ is a finite number. For $-\infty < \sigma < \sigma_0$, curves $\sigma(x, y) = \text{constant}$ intersect a section of the curve $t(x, y) = t^*$ that is outside the branch cut.

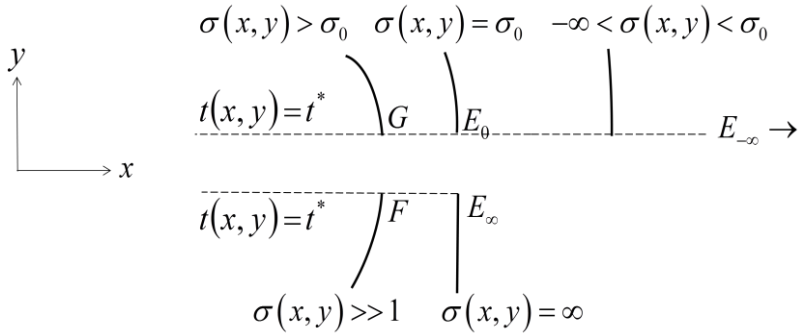

**Fig. 12**   Curves $\sigma = \text{const.}$ for a half-open-ended branch cut (schematic)

Point $E_\infty$, on one side of the end of the branch cut, is intersected by the curve $\sigma(x, y) = \infty$, and point $E_0$, opposite to $E_\infty$, is intersected by the curve $\sigma(x, y) = \sigma_0$. As in the previous case of an open-ended branch cut, the coordinates of $E_\infty$ are given by

$$x(\infty, t^*) = 1,\quad y(\infty, t^*) = 0 \qquad (89)$$

In the limit as $\sigma \to \infty$ along the curve $t(x,y)=t^*$ (i.e. approaching $E_\infty$ on the lower side of the branch cut), it follows from (11) that

$$\frac{\partial x}{\partial \sigma} = -\sum_{n=1}^{\infty} \frac{(-1)^{n-1} \ln n}{n^\sigma} \cos(t^* \ln n) \to 0 , \quad \frac{\partial y}{\partial \sigma} = \sum_{n=1}^{\infty} \frac{(-1)^{n-1} \ln n}{n^\sigma} \sin(t^* \ln n) \to 0 \qquad (90)$$

$$\frac{\partial x}{\partial t} = -\sum_{n=1}^{\infty} \frac{(-1)^{n-1} \ln n}{n^\sigma} \sin(t^* \ln n) \to 0 , \quad \frac{\partial y}{\partial t} = -\sum_{n=1}^{\infty} \frac{(-1)^{n-1} \ln n}{n^\sigma} \cos(t^* \ln n) \to 0 \qquad (91)$$

Now consider points on the upper side of the branch cut in Fig. 12. Since points $E_0$ and $E_\infty$ are infinitesimally close to one another, the coordinates of $E_0$ are also given by (89) as

$$x(\sigma_0, t^*) = 1, \quad y(\sigma_0, t^*) = 0 \qquad (92)$$

The coordinates of any point $G(\sigma, t^*)$ in the neighbourhood of $E_0$ can be obtained by expanding $x(\sigma, t^*)$ and $y(\sigma, t^*)$ about $E_0$ along the curve $t(x,y)=t^*$ in terms of the Taylor's series

$$x(\sigma, t^*) = 1 + (\sigma - \sigma_0)\left(\frac{\partial x}{\partial \sigma}\right)_{\sigma_0, t^*} + O(\sigma - \sigma_0)^2 \ldots \qquad (93)$$

$$y(\sigma, t^*) = (\sigma - \sigma_0)\left(\frac{\partial y}{\partial \sigma}\right)_{\sigma_0, t^*} + O(\sigma - \sigma_0)^2 \ldots \qquad (94)$$

Since $G$ has the same coordinates as a corresponding point $F$ on the opposite side of the branch cut, it follows from (78), (93) and (94) that

$$(\sigma - \sigma_0)\left(\frac{\partial x}{\partial \sigma}\right)_{\sigma_0, t^*} + O(\sigma - \sigma_0)^2 \ldots \approx 0 \qquad (95)$$

$$(\sigma - \sigma_0)\left(\frac{\partial y}{\partial \sigma}\right)_{\sigma_0, t^*} + O(\sigma - \sigma_0)^2 \ldots \approx 0 \qquad (96)$$

For (95) and (96) to be satisfied for *all* values of $\sigma$ in the neighbourhood of $\sigma_0$, we require

$$\frac{\partial x}{\partial \sigma} = \frac{\partial y}{\partial \sigma} = 0 \qquad (97)$$

at $E_0$. It also follows from the Cauchy-Riemann conditions that

$$\frac{\partial y}{\partial t} = \frac{\partial x}{\partial \sigma} = 0 \quad \text{and} \quad \frac{\partial x}{\partial t} = -\frac{\partial y}{\partial \sigma} = 0 \qquad (98)$$

at $E_0$.

(97) and (98) imply that

$$\frac{d\eta}{ds} = \frac{\partial x}{\partial \sigma} + i\frac{\partial y}{\partial \sigma} = \frac{\partial y}{\partial t} - i\frac{\partial x}{\partial t} = 0 \tag{99}$$

at $E_0$.

The inverse derivative $ds/d\eta$ is therefore singular at $E_0$, i.e. at

$$\eta(\sigma_0 + it^*) = x(\sigma_0, t^*) + iy(\sigma_0, t^*) = 1 \tag{100}$$

Using Laurent's theorem, we expand $ds/d\eta$ in a small annular neighbourhood about the singularity at $E_0$, in terms of the series

$$\frac{ds}{d\eta} = \sum_{n=0}^{\infty} c_n (\eta - 1)^n + \frac{d_1}{(\eta - 1)} + \frac{d_2}{(\eta - 1)^2} + \ldots \tag{101}$$

where $c_n$ and $d_n$ are complex constants.

Integrating both sides of (101) with respect to $\eta$, we obtain

$$s(\eta) = \sum_{n=0}^{\infty} \frac{a_n}{n+1}(\eta - 1)^{n+1} + b_1 \ln(\eta - 1) - \frac{b_2}{(\eta - 1)} - \ldots + C \tag{102}$$

where $C$ is a complex constant. (102) implies that $s(\eta)$ is singular in the limit $\eta \to 1$, so that $s(1) = \sigma_0 + it^*$ is not a finite quantity. This contradicts the assumption that $\sigma_0$ and $t^*$ lie in the intervals $-\infty < \sigma_0 < \infty, 0 \leq t^* < \infty$. Hence, the possibility of a half-open-ended branch cut is also rejected.

**Case II(iii)    Closed branch cut(s)**

Finally, consider the closed branch cut, $ABCDA$, shown in Fig. 13, which follows a section of the curve $t(x,y) = t^*$ intersected by curves $\sigma(x,y) =$ constant over the interval $\sigma_C \leq \sigma \leq \sigma_A$

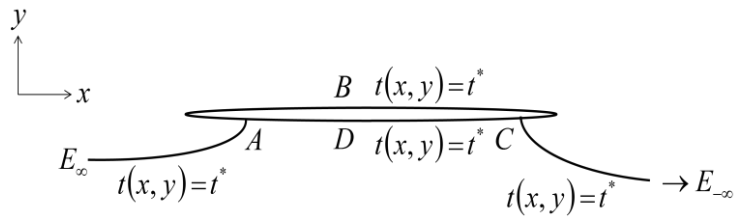

**Fig. 13**        **Closed branch cut**

Curves $\sigma(x,y) = $ constant in the intervals $-\infty < \sigma \leq \sigma_C$ and $\sigma_A \leq \sigma < \infty$ intersect the curve $t(x,y)=t^*$ along the side shoots $CE_{-\infty}$ and $AE_\infty$, respectively, where $E_{-\infty}$ lies on the curve $\sigma(x,y)=-\infty$ and $E_\infty$ lies on the curve $\sigma(x,y)=\infty$. [N.B. The possibility that point $E_\infty$ coincides with point $A$ on the branch cut is not excluded.]

By Theorem 4, $\sigma(x,y)$ decreases monotonically from $E_\infty$ to $E_{-\infty}$ via segment $ADC$, and also decreases monotonically from $E_\infty$ to $E_{-\infty}$ via segment $ABC$. Hence, $\sigma(x,y)$ decreases monotonically from $A$ to $C$ via $D$ and increases monotonically from $C$ to $A$ via $B$. This implies that $\sigma(x,y)$ is non-monotonic around the equipotential following the branch cut, in contradiction to Theorem 4! The possibility of a closed branch cut is therefore ruled out.

The same inconsistency arises if there is more than one closed branch cut, joined by curves $t=t^*$, as shown in Fig. 14.

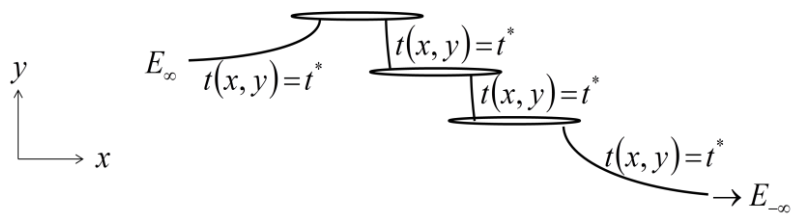

**Fig. 14**     **Multiple closed branch cuts**

## 5     Conclusion

It follows from §4 that none of the possible configurations of curves $\sigma(x,y)=\sigma^*$, $\sigma(x,y)=1-\sigma^*$ and $t(x,y)=t^*$ is compatible with the existence of a zero of $\eta(s)$ in the critical strip for $\sigma^* \neq \tfrac{1}{2}$, but that there no such incompatibility in the case of zeros on the critical line. We conclude that all the zeros of $\eta(s)$ in the critical strip are on the critical line, $\sigma = \tfrac{1}{2}$, and likewise for $\zeta(s)$, consistent with the Riemann hypothesis.

The arguments based on orthogonality (Theorem 2), invariance of the 'flux' integral between stream-lines (Theorem 3) and monotonicity (Theorem 4), clearly apply to other holomorphic functions, but they do not necessarily lead to the same conclusions about the existence of zeros. From §4, the asymptotic behaviour of $\eta(s)$ for $\sigma >> 1$ and $|-\sigma(x,y)| >> 1$ precludes a self-consistent branch cut. However, for holomorphic functions which do not exhibit such asymptotic behaviour, there is no problem in defining a branch cut and zeros on the branch cut are therefore not excluded. An instructive example is given in the Appendix.

## 6     Acknowledgments

The author is indebted to David Ashton, John Norbury and Bruce Pilsworth for their constructive criticism.

**Appendix (see §5)**

We now analyse an example of a holomorphic function which, unlike the Dirichlet eta function, allows a self-consistent open-ended branch cut to be defined and has zeros for $\sigma \neq \frac{1}{2}$.

Consider the function

$$f(s) = s(1-s) = \sigma(1-\sigma) + t^2 + i(1-2\sigma)t \equiv x(\sigma,t) + i\,y(\sigma,t) \tag{A1}$$

where $s = \sigma + it$. $x(\sigma,t)$ and $y(\sigma,t)$ satisfy the Cauchy-Riemann conditions

$$\frac{\partial x}{\partial \sigma} = \frac{\partial y}{\partial t} = 1 - 2\sigma, \qquad \frac{\partial x}{\partial t} = -\frac{\partial y}{\partial \sigma} = 2t \tag{A2}$$

$f(s)$ satisfies the functional equation $f(s) = k(s)f(1-s)$, where $k(s) \equiv 1$, and the symmetry conditions $x(\sigma,t) = x(1-\sigma,-t)$, $y(\sigma,t) = y(1-\sigma,-t)$.

Curves $\sigma(x,y) = \text{constant}$ (stream-lines) are the parabolas

$$y^2 = (1-2\sigma)^2[x - \sigma(1-\sigma)] \tag{A3}$$

and curves $t(x,y) = \text{constant}$ (equipotentials) are the parabolas

$$y^2 = t^2(4t^2 - 4x + 1) \tag{A4}$$

The Riemann surface for $f(s) = s(1-s)$ can be represented on two *x-y* planes, separated by an open-ended branch cut following the degenerate parabola $t(x,0) = 0$ along the *x*-axis over the interval $-\infty < x < \frac{1}{4}$, as shown in Figs. 15(a) and 15(b). Every curve $\sigma(x,y) = \text{constant}$ intersects the branch cut.

The curves $t(x, y) \geq 0$ and $t(x, y) \leq 0$ are situated in the upper and lower planes, respectively. The curves $\sigma(x, y) = \text{constant}$ exist in both planes and cross the branch cut. $\sigma(x,0) = \frac{1}{2}$ is a degenerate parabola which follows the positive x-axis along the strip $\frac{1}{4} \leq x < \infty$.

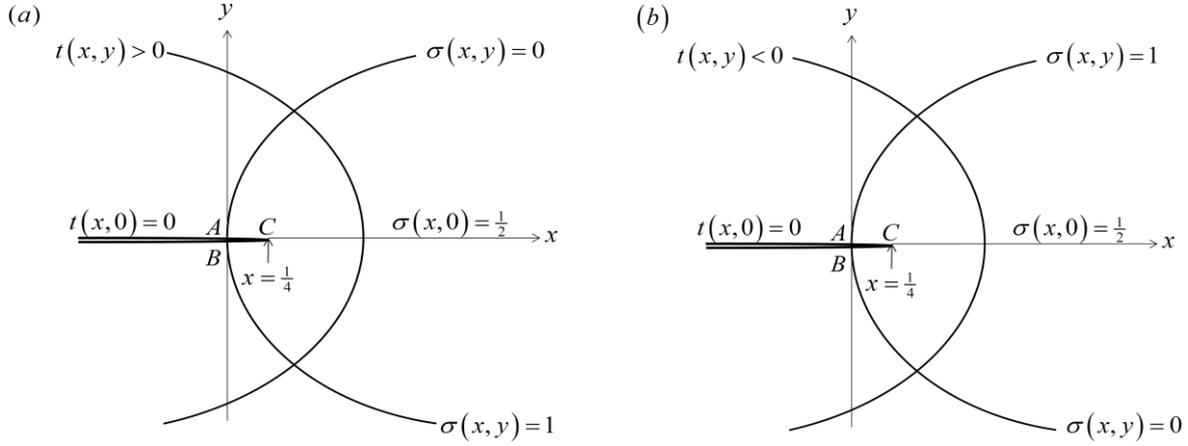

**Figs. 15(a), 15(b)**    Curves $\sigma(x, y) = \text{constant}$ and $t(x, y) = \text{constant}$ in the upper and lower *x-y* planes (schematic)

The zeros of $f(s)$ are situated at $s = 0$ and $s = 1$, at points $A$ and $B$ on the branch cut, at the intersection of the curve $t(x,0) = 0$ with the curves $\sigma(x, y) = 0$ and $\sigma(x, y) = 1$.

The inverse derivative

$$\frac{ds}{df} = \frac{1}{1 - 2s}$$

is singular at $s = \frac{1}{2}$, at the tip of the branch cut ($x = \frac{1}{4}$, $y = 0$).

At any point along the branch cut we have $y = 0$ and $\sigma(x,0) \neq \frac{1}{2}$. (A3) reduces to

$$x = \sigma(1 - \sigma), \quad \text{or} \quad \sigma^2 - \sigma + x = 0 \tag{A5}$$

There are two solutions for $\sigma(x,0)$ for any given value of $x$, given by

$$\sigma(x,0) = \frac{1}{2} \pm \frac{1}{2}\sqrt{1 - 4x} \tag{A6}$$

As $x \to -\infty$, $\sigma(x,0) \to \infty$ on one side of the branch cut and $\sigma(x,0) \to -\infty$ on the opposite side of the branch cut, so there is no inconsistency at the open end of the branch cut in this case.